# *Life without "choice"*

## *The Role of the Axiom of Choice (AC)*

M. A. Sofi

The issue involving the axiom of choice (AC) as a fundamental axiom of mathematics is seen as pesky and uncomfortable within the mathematical community as it's known to lead to such paradoxes as the Banach-Tarski (BT), the existence of a nonmeasurable set of reals or the existence of a Hamel basis in a (non-zero) vector space which, however, is impossible to locate in an uncountable dimensional vector space. The surprising nature of the (BT) Paradox arises from the possibility of cutting up an apple into a finite number of pieces and then juggle them around to get two apples the same size as the original one! This and many more of similar absurd looking consequences of the (AC) had led Roger Penrose to lament that:

***"Most mathematicians would probably regard the (AC) as 'obviously true', while others may regard it as a somewhat questionable assertion which might even be false (and I am myself inclined, to some extent, towards this second viewpoint)".***

Such counterintuitive statements as the existence of non-measurable sets of reals or of the BT-Paradox as consequences of the (AC), would get the sceptics to hope that the world of mathematics would be better off without the axiom of choice being accepted as a natural foundational principle in mathematics. However, as we will learn, that is pure wishful thinking: there arise peskier paradoxes that arise even in the absence of (AC)! (See Section 5).

Now how does one make sense of these counterintuitive statements popping up from all directions: some following as a consequence of a certain axiom (P) while others surface in the absence of (P). That lands us in a typical catch 22 situation!

The main objective of this note is to point out the difficulties that we encounter in an effort to propose, wherever possible, an (AC)-independent proof of certain statements that are known to follow easily from (AC). One such example of the latter part of the assertion is provided by Vitali's famous example of a non-measurable set whose existence is proved as a straightforward consequence of (AC). However, it soon becomes clear how an effort to devise ways to produce such examples while avoiding the use of (AC) gets bogged down in technical difficulties which sometimes look

---

2020 *Mathematics Subject Classification*. Primary 46A22; Secondary 97E60
*Key words and phrases*. Axiom of choice; Lebesgue measure; Invariant measure; Hahn-Banach theorem.

daunting, more so when the idea is to use appropriate consequences of (AC) which are known to be only slightly, but strictly weaker than (AC). The point is that an effort to bypass the use of (AC) often necessitates the need to draw upon methods and tools from many different areas of mathematics to get the desired result. In what follows, this will be shown to be the case in proposing an (AC)-independent proof of the existence of a non-measurable set of reals by invoking the Hahn-Banach theorem which is known to be strictly weaker than (AC) ([4], Theorem 1). For unexplained terms and further insights into certain issues surrounding the Lebesgue measure and the (BT) paradox, we strongly recommend the earlier Resonance articles [5] and [8].

What is the axiom of choice (AC), after all?

In simple terms, the axiom of choice asserts that given an arbitrary family of non-empty sets, there exists a set which consists of exactly one element drawn from each of the given sets.

Slightly more technically, the (AC) stipulates that every given collection of non-empty sets $\mathcal{A} = \{A_\alpha\}_{\alpha\in\Lambda}$ admits a 'choice function', i.e. there exists a function $f: \mathcal{A} \to \coprod_{\alpha\in\Lambda} A_\alpha$ such that $f(A_\alpha) \in A_\alpha$, for all $\alpha \in \Lambda$. An equivalent formulation of (AC) which is widely used across mathematics is provided by the celebrated Zorn's lemma:

Given a partially ordered set (POSET) $(S, \leq)$ such that every chain in it has an upper bound, then S contains a maximal element.

Here, by a *chain* in S, we mean a subset C of S such that any pair of its members are comparable: $for\ all\ a, b \in C$, either $a \leq b$ or $b \leq a$, whereas $a \in S$ is said to be a *maximal* element of S if it is not dominated by any other member of the set: $a \leq b\ for\ some\ b \in S$ implies that $a = b$. Of the many applications of (AC), or of its equivalent variants, the Hahn-Banach theorem of functional analysis comes across as an important example. The (extension form of the) Hahn-Banach theorem guarantees that every bounded linear functional defined on a subspace of a normed linear space admits a bounded linear extension on the whole space. We shall use an equivalent formulation of this theorem in the setting of Boolean algebras as encountered in Theorem 4.1 below

We begin with a discussion of various approaches to a proof of the existence of non-measurable sets in Euclidean spaces that make use of certain (formally) weaker consequences of (AC).

## 1. *Existence of non-measurable sets (AC)*

We begin with

- ***Vitali's proof***

Let $\mathbb{Q} = \{r_n\}$ be an enumeration of rational numbers in [0, 1] and define an equivalence relation '~' on [0, 1] by: x ~ y if x - y ∈ $\mathbb{Q}$. Let $\{M_\alpha\}$ be the collection of equivalence classes induced by ~ and let A be the 'choice' set: $A = \{x_\alpha : x_\alpha \in. M_\alpha\}$. The A is non-measurable. For otherwise, since we have,

$$[0,1] \subseteq \cup_{n=1}^{\infty}(A + r_n) \subseteq [0,2],$$

the translation property of the Lebesgue measure $\mu$ gives

$$1 = \mu([0,1]) \leq \sum_n \mu(\ A + r_n) = \sum_n \mu(\ A) \leq \ \mu([0,2]) = 2,$$

which leads to the absurd conclusion $0 = \mu(A) \neq 0$.

Incidentally, the above argument also shows that there exists no countably additive translation invariant measure $\nu$ on $P(\mathbb{R})$ such that $\nu([0,1]) = 1$.

- ***Via Cauchy's functional equation.***

$$f(x + y) = f(x) + f(y), \text{ for all x, y} \in \mathbb{R}. \qquad (1)$$

This is the so-called Cauchy's functional equation which has been an object of intensive study since it appeared for the first time in a paper of Hamel way back in 1905.

Proposition 1.1: There exist nontrivial solutions to (1).
Proof: Choose a Hamel basis $\{x_\alpha\}_\wedge$ of $\mathbb{R}$ over $\mathbb{Q}$. Then each $x \in \mathbb{R}$ can be written as

$$x = \sum_{\wedge_0} \lambda_\alpha \ x_\alpha$$

where $\wedge_0$ is a finite subset of $\wedge$. For each $\alpha \in \wedge$, define $f_\alpha(x) = \lambda_\alpha$ , $x \in \mathbb{R}$. Then $f_\alpha$ satisfies (1). However, $f_\alpha$ is not continuous, for, $f_\alpha(0) = 0$ , $f_\alpha(x_\alpha) = 1$, and since $f_\alpha$ does not assume irrational values, $f_\alpha$ does not have the intermediate value property. □

Remark 1.2: For f given by the above proposition, it follows that for no c ∈$\mathbb{R}$, it holds that f(x) $= cx, \forall x \in \mathbb{R}$. In fact, it is not difficult to show that this holds precisely when f is continuous. Such solutions shall be referred to as trivial solutions.

We show below that Lebesgue measurable solutions of (1) are trivial. Recall that a function $f: \mathbb{R} \to \mathbb{R}$ is said to be measurable if for each open subset O of $\mathbb{R}$, $f^{-1}(O)$ is measurable.

Proposition 1,3: Given $f: \mathbb{R} \to \mathbb{R}$ satisfying (1) which is Lebesgue measurable, we have $f(x) = cx, x \in \mathbb{R}$, for some c ∈ $\mathbb{R}$.

Proof: The main idea of the proof consists in using Steinhauss theorem (see [1], Theorem 18.13) which asserts that for a measurable set S of reals of positive Lebesgue measure, the difference set S - S contains an open interval about the origin.

Now we can write

$$\mathbb{R} = \bigcup_{n=1}^{\infty} f^{-1}(-n, n).$$

Thus, there exists $\ell \geq 1$ such that m $(f^{-1}(-\ell, \ell)) > 0$, where m denotes the Lebesgue measure. By Steinhauss theorem, there exists an $\alpha > 0$ such that

$$(-\alpha, \alpha) \subset f^{-1}(-\ell, \ell) - f^{-1}(-\ell, \ell).$$

Since f is additive, we have

$$f(-\alpha, \alpha) \subset (-2\ell, 2\ell).$$

This yields that f is continuous. For, given $\varepsilon > 0$, choose N $\geq 1$ such that $\frac{2\ell}{N} < \varepsilon$. Now, for $|x| < \delta = \frac{\alpha}{N}$, we have

$$|f(Nx)| = N|f(x)| < 2\,\ell.$$

This gives

$$|f(x)| < \frac{2\ell}{N} < \varepsilon$$

which means that f is continuous at the origin and hence continuous throughout. In other words, there exists c ∈ ℝ such that f(x) = cx, for all x∈ ℝ.

□

The above result yields that the nontrivial solution f of (1) as guaranteed by Proposition 1.1 combined with Remark 1.2 is non-measurable (in the Lebesgue sense). Hence there exists an open subset $O$ of ℝ such that $f^{-1}(O)$ is not measurable in R. However, it is possible to show that there exist a $\sigma -$ $a$lgebra that properly contains the Lebesgue measurable sets with respect to which these nontrivial solutions are measurable.

## *2. Towards non-measurability by other means*

Apart from the use of (AC) in the proof of the existence of non-measurable sets as described in the previous section, we describe below yet another method involving the use of certain related, but well-known paradoxes to provide further examples of such sets. To this end, we shall see how the proof of the Hausdorff paradox can be modified to avoid the use of (AC) in the argument and thus produce non-measurable sets without the use of (AC). Before we do that (Section 4), it is useful to illustrate the idea of a paradox and illustrate it with the help of some examples.

We shall use the following notation
$\mathbb{R}^n$: Euclidean space of dimension n.
$S^n$: Unit sphere in $\mathbb{R}^n$.
$B_n$: Closed unit ball in $\mathbb{R}^n$.
$O_n$: The group of (orthogonal) n$\times n$ matrices $A$ such that $AA^t = A^tA = I$ where $A^t$ is the transpose of $A$ and $I$ is the identity matrix.
$SO_n$: Special orthogonal group: The subgroup of $O_n$ with determinant equal to1.
$isom(\mathbb{R}^n)$: The group of bijective isometries on $\mathbb{R}^n$.

Definition 2.1: A group action induced by a group G on a set X is given by a mapping

$$G \times X \to X, (g, x) \to g(x) \in X.$$

such that for each $g \in G$, the map $x \to g(x)$ acts as a bijection on X. More precisely, a group action (G, X) is given by a group homomorphism between the group G and the group of (bijective) transformations on X such that $(g.h)x = g.(hx)$ and $e(x) = x$ for all $g, h \in G, x \in X$. The action is said to be *free* if $(g, x) \to x$ implies that $g = e$, the identity of G.

Definition 2.2: Given a group G and a G-action on a set X, we say a subset E ⊆ X is G-*paradoxical* (or paradoxical with respect to G) if there exist pairwise disjoint subsets $A_1, A_2, \dots, A_m, B_1, B_2, \dots, B_n$ of E and $g_1, g_2, \dots, g_m, h_1, h_2, \dots, h_n \in G$ such that E $= \bigcup_{i=1}^{m} g_iA_i = \bigcup_{i=1}^{n} h_iB_i$. Further, a group is said to be paradoxical if it is paradoxical under the action of (left) multiplication by the elements of the group itself.

Definition 2.3: Under the conditions of the previous definition, given subsets E, F ⊆ X, we say E and F are G-*equidecomposable* if there exist pairwise disjoint subsets $A_1, A_2, \dots, A_m$ *of E* and $g_1, g_2, \dots, g_m \in G$ such that E $= \bigcup_{i=1}^{n} A_i$ and F $= \bigcup_{i=1}^{n} g_iA_i$.

Theorem 2.4: There exists a set in the plane which is paradoxical under the isometry group.

Proof: We note that given distinct polynomials P, Q with rational coefficients and a transcendental number t, then P(t) $\neq Q(t)$. Consider the set

$$\mathrm{E} = \{P(t): P \in \mathbb{Z}_+[x]\}.$$

Then E is paradoxical. Indeed, let A = $\{P(t): P \in \mathbb{Z}_+[x], P(0) = 0\}$ and B = $\{P(t): P \in \mathbb{Z}_+[x], P(0) \neq 0\}$. Clearly, A and B are disjoint and $E = A \cup B$. Finally, taking $t = e^i$ and defining the isometries by $g(a) = t^{-1}a$, and $h(b) = b - 1$, we see that $g(A) = E$ and $h(B) = E$. Indeed, let $P \in E$, then $tP \in A$, $P + 1 \in B$ and, therefore, $g(A) = E, h(B) = E$. ◻

Theorem 2.4 is subsumed within the following more general assertion:

Theorem 2.5: Let G be a group acting on a set X such that G contains $\rho, \tau$ such that for some x in X, any two words in $\rho, \tau$ and beginning with $\rho, \tau$, respectively, differ when applied to x. Then X contains a nonempty G-paradoxical set. (The desired paradoxical set is obtained by taking the G-orbit at x). (See [10], Theorem 1.8 for a proof).

Comments 2.6: An important step towards the search for non-measurable sets is given by the following facts (which are of independent interest) and which are used in the 'choice-based' proof of the Hausdorff paradox that $S^2$ is $SO_3$ −paradoxical.

a. Free group on two generators $F_2$.

A *free* group generated by a set S is the set of all *words* in the alphabet $\mathcal{A} = \{s, s^{-1}: s \in S\}$, where a *word* is a string of symbols $x_1, x_2, \ldots, x_n$ with each "*letter*" $x_i \in \mathcal{A}$. A *reduced word* is one where $s$ *and* $s^{-1}$ do not stand next to each other. The free group G with S as its *generating set* is the set of all reduced words in the alphabet $\mathcal{A} = \{s, s^{-1}: s \in S\}$ under *concatenation*: given reduced words a $= x_1x_2 \ldots x_n$ and b $= y_1y_2 \ldots y_n$, we define a$* b = x_1x_2 \ldots x_n y_1 y_2 \ldots y_n$, with the *empty word* – containing no symbols – representing the identity of G. The *rank* of G is the cardinality of the (minimal) generating set S. The free group of rank 2 shall be denoted by $F_2$.

b. $F_2$ is paradoxical.

Proof: Let a, b be the generators of $F_2$. For $x \in \{a, b, a^{-1}, b^{-1}\}$, let $W(x) = \{w \in F_2: w \textit{ starts with } x\}$. Clearly, we have

$$F_2 = \{e_{F_2}\} \cup W(a) \cup W(b) \cup W(a^{-1}) \cup W(b^{-1}).$$

It follows that

$$F_2 = W(a) \cup aW(a^{-1}) \text{ and } F_2 = W(b) \cup bW(b^{-1}).$$

This is so because if $h \in F_2 \backslash W(a)$, then $a^{-1}h \in W(a^{-1})$. This gives $h = a(a^{-1}h) \in aW(a^{-1})$. □

$c. SO_3$ contains $F_2$ as a subgroup. In other words, $SO_3$ contains a subgroup which is isomorphic with $F_2$.

As free generators, one may consider the matrices

$$A = \frac{1}{7}\begin{pmatrix} 6 & 2 & 3 \\ 2 & 3 & -6 \\ -3 & 6 & 2 \end{pmatrix}, B = \frac{1}{7}\begin{pmatrix} 2 & -6 & 3 \\ 6 & 3 & 2 \\ -3 & 2 & 6 \end{pmatrix}$$

and show that the subgroup generated by A and B (as a free group of rank 2) is isomorphic to $F_2$. See [10], Theorem 2.1.

Proposition 2.7 (AC) ([10], Theorem 1.10): If a group acts freely on a set X such that G is paradoxical, then X is G-paradoxical. In particular, if a group G contains a paradoxical subgroup, then G is also paradoxical.

Corollary 2.8: $SO_3$ is paradoxical.

The preceding statement follows by combining (b) and (c) with Proposition 2.7.

To deduce the Hausdorff Paradox from the previous assertions, let us note that if $SO_3$ were to act freely on $S^2$, then an application of Proposition 2.7 in combination with Corollary 2.8 would lead us to conclude that $S^2$ is indeed $SO_3 -$ paradoxical. Unfortunately, the action is obviously not free. Note that each element M of $SO_3$ - being an orthogonal matrix - has 1 as an eigenvalue, so each element of the 1-dimensional eigenspace of M corresponding to the eigenvalue 1 remains fixed under the action of M. In other words, the action of M on $\mathbb{R}^3$ amounts to a rotation of the unit sphere $S^2$ that leaves exactly two points on $S^2$ fixed which lie on the axis of the given rotation. Thus Proposition 2.7 cannot be applied, at least in this form, to derive the desired conclusion.

To circumvent this obstacle, consider the countable set C of all such fixed points to get a free action by the subgroup $F_2$ of $SO_3$ on $S^2\backslash C$. Now it is legitimate to apply Proposition 2.7 to $S^2\backslash C$ to conclude, courtesy Corollary 2.8, that $S^2\backslash C$ is $F_2 -$ paradoxical. The final push is provided by the $SO_3 -$ equidecomposability of $S^2\backslash C$ and $S^2$ and the fact that paradoxical character is preserved under equidecomposability to conclude that $S^2$ is $SO_3 -$ paradoxical. As we shall see later (Corollary 4.4), this leads to the existence of a nonmeasurable subset of $S^2$.

In an effort to produce non-measurable sets without (AC), we follow the (AC)-independent approach as proposed by Foreman and Wehrung [2] where we make use of a suitable version of the Hahn Banach theorem which is known to be (strictly) weaker than (AC).

## *3. Banach meets Lebesgue*
## *Hahn-Banach theorem and the existence of non-measurable sets.*

We begin with the following definition.

Definition 3.1: A discrete group (semigroup-a set S equipped with an associative binary operation) S is said to be ***(right) amenable*** if there exists a (right) ***invariant mean*** on $B(S)$, i.e., a (norm-one) linear functional F on $B(S)$ – the Banach space of all bounded functions on S with the uniform norm - such that F(e) =1 and F(f) = F($f^s$) where $f^s(t) = f(ts), t, s \in S$. Similarly for the left invariant mean. Further, S is said to be ***amenable*** if it is both right and left amenable.

Remark 3.2: It is easily seen that the above condition leads to the existence of a *Banach measure* on a given amenable semigroup S - a translation invariant, finitely additive probability measure defined on all subsets of S via: $\nu(A) = F(\chi_A)$. In fact, converse is also true (check!).

Example 3.3: Finite groups are amenable: take $\mu(A) = {|A|}/{|G|}$ where G is a finite group and A⊆G.

Example 3.4: The group $\mathbb{Z}$ is amenable.

Proof: Let $\mu_n(A) = \frac{|A\cap\{0,1,\ ...\ ,n-1|}{n}$. Note that $\mu_n$ is a finitely-additive measure on $\mathbb{Z}$, $0 \leq \mu_n(A) \leq 1$ and $|\mu_n(A+1) - \mu_n\ (A)| \leq {2}/{n}$ for all A⊆$\mathbb{Z}$ where we use the notation: A + 1:= {x + 1: x ∈ A}. Thus $\mu_n \in [0,1]^{P(\mathbb{Z})}$ which is compact by Tychonoff's theorem - as the product of $|P(\mathbb{Z})|$ many copies of the compact interval [0, 1]- so there exists a convergent subsequence $\{\mu_{n_k}\}$ of $\{\mu_n\}$: $\mu_{n_k} \to \mu$, say. It is easily checked that $\mu$ has all the desired properties. □

**Notation**: Given a group G acting on a set X, we shall use the symbol (G, X) whereas the notation IM (G, X) shall be used to mean the following statement:

(*) There exists a *G-invariant* finitely additive measure $\mu$ on P(X):
( $\mu(A) = \mu(\mathrm{g}^{-1}A)$, for all $A \in P(X)$ and for all $g \in G$) such that $\mu(X) = 1$.

A simple example of the above property is provided by amenable groups as already noted in Remark 3.2 above. On the other hand, the existence of invariant Borel measures on a large class of topological spaces is guaranteed by the following statement.

Proposition 3.5: Each compact Hausdorff space X admits a regular Borel measure which is invariant under a given continuous mapping on X.

The proof follows by a delightful application of the existence of Banach limts – a special case of the invariant mean - on $\ell_\infty$ (of bounded scalar sequences) combined with the Reisz representation theorem on the dual of C(X), the Banach space of scalar valued continuous functions on a compact Hausdorff space X.

Definition 3.6: A *Banach limit* is a positive shift invariant continuous linear functional L on the space $\ell_\infty$ such that L(e) = 1 (e = (1,1,1, . , . . ,1, , . ,).

Banach limits always exist. For a proof based on the Hahn-Banach theorem, see [6]. A straightforward consequence of the existence of Banach limits is that the semigroup $\mathbb{N}$ is amenable. A suitable modification of the proof for the existence of Banach limits on the space $B(\mathbb{R}) - the\ space\ of\ bounded\ functions\ on\ \mathbb{R} -$ yields that $\mathbb{R}$ is an amenable group: the measure $\mu$ witnessing the amenability may be defined by the formula as indicated in Remark 3.2 above.

On the other hand, the Reisz representation theorem provides a complete description of all the bounded linear functionals on C(X):

A linear functional F on C(X) is continuous if and only if there exists a unique regular Borel measure $\mu$ on X such that $F(f) = \int f d\mu$, for each $f \in C(X)$.

For an accessible proof, we recommend [7].

Proof of Proposition 3.5: Let $T: X \to X$ be a given continuous mapping and let $T^n$ denote the nth iterate of T and consider the group G generated by T. Our proof yields the stronger conclusion that the measure is in fact invariant with respect to (the semigroup) of all powers $T^n$, $n \geq 1$. Choose a Banach limit $\Lambda$ on $\ell_\infty$ and fix some $x \in X$. Consider the linear functional F on C(X) given by

$$F(f) = \Lambda\big(f(x), f\big(T(x)\big), f\big(T^n(x)\big), \dots\big).$$

An application of the Reisz Representation Theorem yields a unique regular Borel measure $\mu$ on X such that $F(f) = \int f d\mu$ for each $f \in C(X)$. Using the (shift) invariance of $\Lambda$ gives

$$\Lambda\big(f(x), f\big(T(x)\big), f\big(T^3(x)\big), \dots\big) = \Lambda\big(f\big(T(x)\big), f\big(T^2(x)\big), f\big(T^3(x)\big), \dots\big),$$

or equivalently, $\int f d\mu = \int f o T d\mu$ for each $f \in C(X)$. Finally, invoking the change of variables formula from calculus, we see that $\int f d\mu = \int f d\mu T^{-1}$ holds for each $f \in C(X)$. Finally, the uniqueness of the measure $\mu$ gives $\mu = \mu T^{-1}$, yielding thereby that $\mu$ is T-invariant. A repeated application of shift invariance shows that in fact, $\mu$ is invariant under $T^n$ for all $n > 1$. □

Comment: The same argument as used above can be used to show the invariance of $\mu$ with respect to the group of all the powers $T^n, n \in \mathbb{Z}$ provided T is assumed to be a homeomorphism.

Proposition 3.7: G, Amenable ⇒ IM (G, X) ⇒ X is G-nonparadoxical.

Proof: The first implication is a direct consequence of the definition: given a finitely additive invariant probability measure $\mu$ on G and a fixed $x \in X$, then $\nu$ given by $\nu(E) = \mu\{g \in G: g(x) \in E\}$ defines the desired measure on X, thus verifying IM (G, X).

To prove the second implication, assume, on the contrary, that X is G-paradoxical and choose pairwise disjoint subsets $A_1, A_2, \dots, A_m, B_1, B_2, \dots, B_n$ of X and $g_1, g_2, \dots, g_m, h_1, h_2 \dots h_n \in G$ such that $E = \bigcup_{i=1}^m g_i A_i = \bigcup_{i=1}^n h_i B_i$. Let $\nu$ be a measure witnessing IM (G, X). We have

$$\nu(X) \geq \sum_{i=1}^m \nu(A_i) + \sum_{1=1}^n \nu(B_i)$$

$$= \sum_{i=1}^m \nu(g_i A_i) + \sum_{1=1}^n \nu(h_i B_i)$$

$$\geq \nu\left(\bigcup_{i=1}^m g_i A_i\right) + \nu\left(\bigcup_{i=1}^n h_i B_i\right)$$

$$= \nu(X) + \nu(X)$$

yielding $\nu(X) = 0$, a contradiction. □

Combined with 2.6(b), Proposition 3.7 gives

Corollary 3.8: $F_2$ is not amenable.

## 4. Non-measurable sets without (AC)

*Bringing Hahn-Banach theorem to the table*

In what follows, we tailor the familiar traditional approach to a proof of the Hausdorff paradox to create a setup where it becomes possible to avoid the use of (AC) via an application of (an appropriate version of) the Hahn Banach theorem to derive certain consequences on a group that acts freely on a universal measure space. Here, we use an analogue of the Hahn Banach theorem for measures defined on a Boolean algebra stated below.

Theorem 4.1 (ZF + Hahn Banach [10]): Every Boolean algebra admits a [0, 1]-valued finitely additive probability measure.

We recall that a Boolean algebra is a quadruplet $(\mathcal{A}, \vee, \wedge, \sim)$ consisting of a set $\mathcal{A}$ equipped with three operations $\vee, \wedge$ and $\sim$ where $\sim$ is a unary operation while $\vee$ and $\wedge$ are binary operations on $\mathcal{A}$ satisfying the following axioms for all $a, b \in \mathcal{A}$:

(i) Both $\vee$ and $\wedge$ are commutative, associative and distributive.
(ii) $(a \wedge \mathrm{b}) \vee \mathrm{b} = (a \vee \mathrm{b}) \wedge \mathrm{b} = \mathrm{b}$.
(iii) $(a \wedge a^{\sim}) \vee \mathrm{b} = (a \vee a^{\sim}) \wedge \mathrm{b} = \mathrm{b}$.

Here, $a^{\sim}$ is called the complement a in $\mathcal{A}$. By virtue of (iii), it follows that both $a \wedge a^{\sim}$ and $a \vee a^{\sim}$ define elements of $\mathcal{A}$ which are independent of a and are denoted by 0 and1, respectively: $a \wedge a^{\sim}=0$, $a \vee a^{\sim} = 1$. Further, we say that $a, b \in \mathcal{A}$ are disjoint if $a \wedge \mathrm{b} = 0$.

A standard example of a Boolean algebra is provided by the power set P(X) of any given set X where the Boolean operations are provided by the union, intersection and the complement of sets.

A subalgebra of a Boolean algebra $\mathcal{A}$ is a subset $\mathcal{A}_0$ which is closed under $\vee, \wedge$ and $\sim$. Finally, a measure on $\mathcal{A}$ is a function $\mu: \mathcal{A} \to [0, \infty]$ such that $\mu(0) = 0$ and $\mu(a \vee \mathrm{b}) = \mu(a) + \mu(b)$ if $a \wedge \mathrm{b} = 0$.

Given a group G of automorphisms (maps that preserve the three Boolean operations) on $\mathcal{A}$, the measure $\mu$ is said to be G-invariant if $\mu(a) = \mu(ga)$ for all $a \in \mathcal{A}$, $\mathrm{g} \in G$. By a universal measure space we shall mean the Boolean algebra P(X) on a given set X equipped with a (finitely additive) measure and will be denoted as (X, $\mu$).

In the next theorem, we shall make use of the definition of the integral of a function on a measure space (X, $\Sigma. \mu$) where $\mu$ is a finitely additive measure on $\Sigma$ which is analogous to the standard definition of the Lebesgue integral. Accordingly, we define the integral of a simple function $f$ on X and note that the integral $\int f(x)d\mu$ defines a linear functional on the space $S(\Sigma)$ of simple functions, so it extends (uniquely) to a linear functional on the uniform closure $B(\Sigma)$ of $S(\Sigma)$ (with respect to the sup-norm). Thus, the integral $\int f(x)d\mu$ makes sense for all functions $f \in B(\Sigma)$. It turns out that, as in the case of a Lebesgue integral on a closed interval in $\mathbb{R}$, $f: X \to \mathbb{R}$ is integrable in the above sense if it is bounded and measurable.

Theorem 4.2: Assume that a group G acts freely on some universal measure space (X, $\mu$). Then G is amenable.

Proof: Let $X/_G$ denote the set of all orbits of X induced by G. By Theorem 4.1, there exists a universal probability measure $\mu_x$ on each orbit [x] $\in X/_G$. Thus, $\mu_x$ is a measure on $P([\text{x}])$. For each $A \subset G$, define a function $f_A: X \to [0,1]$ given by $f_A(x) = \mu_x(Ax)$ and use it to define $\sigma: P(G) \to [0,1]$ by $\sigma(A) = \int f_A(x)d\mu$. To show that $\sigma$ is an invariant finite additive measure, let A and B be disjoint subsets of G. Because G acts freely, it follows that $f_C = f_A + f_B$ where $C = A \cup B$. This gives $\sigma(A \cup B.) = \sigma(A) + \sigma(B)$. Further, if $A, B \subset G$, $g \in G$ and B = Ag, then for $x \in X$, we have

$$f_B(x) = \mu_x(Bx) = \mu_x(Agx) = \mu_x(A(gx)) = \mu_{gx}(A(gx)) = f_A(gx).$$

Combining the above with the invariance of $\mu$, we get

$$\sigma(B) = \int f_B(x)d\mu = \int f_A(gx)d\mu = \int f_A(x)d\mu g^{-1} = \sigma(A). \qquad \square$$

Theorem 4.3: Assume that $IM(SO_3, S^2)$ holds. Then there exists a free action of $F_2$ on some universal measure space.

Proof: Under the given hypothesis, there exists a $SO_3$ −invariant finitely additive probability measure $\mu$ on $S^2$. Observe that each non-identity rotation r in $SO_3$ has an axis of rotation which intersects the sphere $S^2$ exactly in two points. Further, these are the only points which remain fixed under the rotation r. Let C be the set of all such points on the sphere which remain fixed corresponding to rotations in the (countable free) subgroup $F_2$ of $SO_3$. Because $F_2$ is countable, it follows that C is countable too, and so we can choose a line L which bypasses the set C. It is clear that $F_2$ acts freely on $S^2 \backslash C$.

To complete the proof, it suffices to show that $\mu(C) = 0$. That would ensure that $\mu(S^2 \backslash C) = 1$ and then we can choose our universal measure space to be $(S^2 \backslash C, \mu)$.

Claim: There exists $g \in SO_3$ such that $g^i(C) \cap g^j(C) = \phi, \forall\, i \neq j$.
It suffices to show that $g^i(C) \cap C = \phi, \forall\, i > 0$, for if $i > j$, then we have

$$g^i(C) \cap g^j(C) = g^j(C \cap g^{i-j}(C)) = g^j(\phi) = \phi.$$

To this end, we show that the set S of $g \in SO_3$ such that $g^i(C) \cap C \neq \phi$ for some $i > 0$ is at most countable. The existence of the desired g would follow from the uncountability of $SO_3$.

For $P \in C$, let A(P) be the set of all angles $\theta$ such that P is mapped onto other points of C by some rotation around the line L described above through an angle $\theta$. Because C is countable, it follows that A(P) is countable, and so is the union $A = \{A(P): P \in C\}$. It follows that if g is a rotation around L through one of the uncountably many angles not in A(P), we have $g^i(C) \cap C = \phi, \forall\, i > 0$. Finally, assume that $\mu(C) > 0$ and choose $n \geq 0$ such that $n\, \mu(C) > 1$. By the $SO_3 -$ invariance of $\mu$, this gives

$$\mu(C) + \mu\big(g(C)\big) + \mu(g^2(C)) + \ldots + \mu(g^n(C) > 1 = \mu\,(S^2),$$

a contradiction. □

We now prove our promised theorem which yields a delightful consequence involving free action of a group on a universal measure space.

Finally, Theorem 4.2 and 4.3 combined with Corollary 3.8 yields

Corollary 4.4: $(IM(SO_3,\ S^2)$ does not hold. In particular, there exist non-measurable subsets of $S^2$ – and hence of $\mathbb{R}^3$.

To see why the latter statement holds, recall that the Lebesgue measure is rotation invariant as a countably (hence finitely) additive measure on the $\sigma -$algebra $M$ of all measurable subsets of $S^2$. Combined with the previous corollary, it follows that the Lebesgue measure on $S^2$ would fail to be a finitely additive $SO_3 -$ invariant measure on $P(S^2)$, thereby yielding that $M \neq P(S^2)$ and, thus, a subset of $S^2$ which is not (Lebesgue) measurable. □

## *5: (AC) is not a bad guy, after all!*

In Section 1, we saw how a straightforward application of the axiom of choice led to the existence of a non-measurable set of real numbers. In the preceding paragraphs, we noted how an attempt to look for such examples in the absence of (AC) was fraught as it was seen to cost considerable effort in drawing the same conclusion by invoking methods which were required to be (AC)-independent. However, as it turned out, that was made possible, but by appealing to a heavy arsenal of sophisticated technical machinery borrowed from certain areas of advanced mathematics. That makes a strong case for avoiding to be squeamish about the use of (AC) in mathematics, a claim that is further buttressed by a slew of paradoxes that arise in mathematics in the absence of the (AC). Below we provide one example of that phenomenon. For details, see [9].

Theorem 5.1(ZF): $|\mathbb{R}| < \left|{}^{\mathbb{R}}/_{\mathbb{Q}}\right|$.

An interpretation of the above is provided by the following commonplace illustration:

"In a sports tournament that features a few teams participating in the event, it's possible that the number of teams far exceeds the total number of players who are representing their teams"! Put differently, the theorem amounts to saying that a country could have more inhabited provinces than it has people!

*5.2: Banach-Tarski Paradox (ZF) and non-measurable sets*

The Banach-Tarski paradox in $\mathbb{R}^3$ - and in higher dimensions - that shows the non-existence of finitely additive $\text{isom}(\mathbb{R}^3) -$ invariant measures on $\mathbb{R}^3$- quickly leads to the existence of non-measurable sets which, as we have seen in the preceding sections, rely on (AC), considering that the traditional proofs of the (BT) paradox are heavily (AC)-dependent. As if on a cue, soon after Foreman and Wehrung [2] appeared, Pawlowski [3] showed the following important result towards a proof of the (BT) paradox using the Hahn-Banach theorem while avoiding the (AC).

Theorem 5.3 (ZF) [3]: If the free group F on two generators acts freely on a set X, then X is F-paradoxical.

The Banach-Tarski paradox as shown above (without (AC)) and combining it with Proposition 3.7 and the argument as used in the proof of Corollary 4.4 yields the existence of a non-measurable subset of $B_3$.

Acknowledgements: The author wishes to thank the referee for carefully reading through the manuscript and for making suggesting changes that have helped present a (maximally possible) self-contained version of the material.

Department of Mathematics, Kashmir University, Srinagar, India.
Email address: aminsofi@gmail.com